\documentclass[12pt,amstex]{amsart}

\usepackage{mathptmx}
\usepackage{mathrsfs}

\usepackage{verbatim}
\usepackage{url}
\usepackage[all]{xy}
\usepackage{color}

\usepackage[colorlinks=true,citecolor=blue]{hyperref}

\usepackage{stmaryrd}
\usepackage{epsfig}
\usepackage{amsmath}
\usepackage{amssymb}

\usepackage{amscd}
\usepackage{graphicx}
\usepackage{pstricks}

\theoremstyle{plain}
\newtheorem{thm}{Theorem}[section]

\newtheorem{cor}[thm]{Corollary}
\theoremstyle{definition}

\theoremstyle{remark}

\def \N {\mathbb N}

\def \Z {\mathbb Z}

\def \O {\mathcal{O}}

\def \ep {\epsilon}
\def \d {\delta}

\def \w {\omega}

\def \lra {\longrightarrow}

\begin{document}

\title{Reducing the Sarnak Conjecture to Toeplitz systems}

\author{Wen Huang}

\author{Zhengxing Lian}
\author{Song Shao}
\author{Xiangdong Ye}

\address{Wu Wen-Tsun Key Laboratory of Mathematics, USTC, Chinese Academy of Sciences,
Department of Mathematics, University of Science and Technology of China,
Hefei, Anhui, 230026, P.R. China and Institute of Mathematics, Polish Academy of Sciences.}

\email{wenh@mail.ustc.edu.cn}
\email{lianzx@mail.ustc.edu.cn}
\email{songshao@ustc.edu.cn}
\email{yexd@ustc.edu.cn}

\subjclass[2010]{Primary: 37B05; 54H20}

\thanks{This research is supported by NNSF of China (11971455, 11571335,  11431012). The second author was partially supported by the National Science Centre (Poland) grant 2016/22/E/ST1/00448.}

\date{}

\begin{abstract}
In this paper, we show that for any sequence ${\bf a}=(a_n)_{n\in \Z}\in \{1,\ldots,k\}^\mathbb{Z}$ and any $\epsilon>0$, there exists a Toeplitz sequence ${\bf b}=(b_n)_{n\in \Z}\in \{1,\ldots,k\}^\mathbb{Z}$ such that the entropy $h({\bf b})\leq 2 h({\bf a})$ and $\lim_{N\to\infty}\frac{1}{2N+1}\sum_{n=-N}^N|a_n-b_n|<\epsilon$.
As an application of this result, we reduce Sarnak Conjecture to Toeplitz systems, that is, if the M\"{o}bius function is disjoint from any Toeplitz sequence with zero entropy, then the Sarnak conjecture holds.
\end{abstract}
\maketitle

\section{Introduction and Preliminaries}

\subsection{Topological dynamical systems}

A pair $(X,T)$ is called a {\em topological dynamical system} (t.d.s. for short) if $X$ is a compact metric space and $T : X \rightarrow X$ is a homeomorphism.
A topological system $(X, T)$ is {\em transitive} if there exists some point
$x\in X$ whose orbit $\O(x,T)=\{T^nx: n\in \Z\}$ is dense in $X$ and
we call such a point a {\em transitive point}. The system is {\em
minimal} if the orbit of any point is dense in $X$.
A point $x\in X$ is called a {\em minimal point} if
$(\overline{\O(x,T)}, T)$ is minimal. A well-known fact is that
$x\in X$ is a minimal point if and only if for any neighborhood $U$ of $x$, the returning time set
$$\mathbf{N}(x,U)=\{n\in \Z:\ T^nx\in U\}$$
is syndetic. i.e. there is $l>0$ such that for each $m\in \Z$, $[m,m+l]\cap N(x,U)\neq \emptyset$, where $[m,m+l]=\{m,m+1,\ldots,m+l\}$ is the interval of $\Z$.

A {\it factor map} $\pi: X\rightarrow Y$ between the t.d.s. $(X,T)$
and $(Y,S)$ is a continuous onto map which intertwines the
actions; one says that $(Y,S)$ is a {\it factor} of $(X,T)$ and
that $(X,T)$ is an {\it extension} of $(Y,S)$. One says that:
$\pi$ is an {\it almost one to one} extension  if there
exists a dense $G_\d$ set $X_0\subseteq X$ such that
$\pi^{-1}(\{\pi(x)\})=\{x\}$ for any $x\in X_0$.

\subsection{Sarnak conjecture}
The M\"{o}bius function $\mu: \N\rightarrow \{-1,0,1\}$ is defined by
$\mu(1)=1$ and
\begin{equation}\label{M-function}
  \mu(n)=\left\{
           \begin{array}{ll}
             (-1)^k & \hbox{if $n$ is a product of $k$ distinct primes;} \\
             0 & \hbox{otherwise.}
           \end{array}
         \right.
\end{equation}

We say a sequence $\xi$ is {\em realized} in $(X,T)$ if there is an $f\in C(X)$ and an $x\in X$
such that $\xi(n) = f(T^nx)$ for any $n\in\N$. A sequence $\xi$ is called {\em deterministic} if it is realized in a
system with zero topological entropy. Here is the conjecture by Sarnak \cite{Sar}:

\medskip

\noindent {\bf Sarnak Conjecture:}\ {\em
The M\"{o}bius function $\mu$ is linearly disjoint from any deterministic sequence $\xi$. That is,
\begin{equation*}\label{Sarnak}
  \lim_{N\rightarrow \infty}\frac{1}{N}\sum_{n=1}^N\mu(n)\xi(n)=0.
\end{equation*}
}

\subsection{Symbolic dynamics}

Let $S$ be a finite alphabet with $k$ symbols, $k \ge 2$. We assume that $S=\{1,2,\cdots,k\}$. Let $\Sigma_k=S^{\Z}$ be the
set of all sequences ${\bf x}=\ldots x_{-1}x_0x_1 \ldots=(x_i)_{i\in \Z}$, $x_i \in
S$, $i \in \Z$, with the product topology. A metric compatible is
given by $d({\bf x},{\bf y})=\frac{1}{1+k}$, where $k=\min \{|n|:x_n \not= y_n
\}$, ${\bf x},{\bf y} \in \Sigma_k$. The shift map $\sigma: \Sigma_k \lra \Sigma_k$
is defined by $(\sigma {\bf x})_n = x_{n+1}$ for all $n \in \Z$. The
pair $(\Sigma_k,\sigma)$ is called a {\em shift dynamical system}. Any subsystem of $(\Sigma_k, \sigma)$
is called a {\em subshift system}.

Each element of $S^{\ast}= \bigcup_{k \ge 1} S^k$ is called {\em a word} or {\em a block} (over $S$).
We use $|A|=n$ to denote the length of $A$ if $A=a_1\ldots a_n$.
If $\omega=(\ldots \omega_{-1} \omega_0 \omega_1 \ldots) \in \Omega$ and $a \le b \in \Z$, then
$\omega[a,b]=\omega_{a} \omega_{a+1} \ldots \omega_{b}$ is a $(b-a+1)$-word occurring in
$\omega$ starting at place $a$ and ending at place $b$.
Similarly we define $A[a,b]$ when $A$ is a word. A finite-length word $A$ {\em appears} in the word $B$ if there are some $a\le b$ such that
$B[a,b]=A$.

Let $(X,\sigma)$ be a subshift system. The collection of all $n$-words of $X$ is denoted by $B_n(X)$. Then the {\em topological entropy} of $(X,\sigma)$ is defined by
$$h(X,\sigma)=\lim_{n\to\infty}\frac{\log \# B_n(X)}{n},$$
where $\# $ means cardinality.

Let ${\bf \omega}\in \Sigma_k$. Denote the orbit closure of ${\bf \omega}$ by $X_{\bf \omega}$. We call $h(X_{\bf \omega}, \sigma)$ the {\em entropy} of the sequence ${\bf \omega}$, and also denoted it by $h({\bf \omega})$.

\subsection{Toeplitz systems}

Fixing a sequence of positive integers $(n_i)_{i\in \Z}$ with with $n_i\ge 2$, let $X=\prod_{i=1}^\infty\{0,1,\ldots, n_i-1\}$ with the discrete topology on each coordinate and the product topology on X. Let $T$ be the transformation that is addition by $(1,0,0,\ldots)$ with carrying to the right. Thus $T(x_1,x_2,\ldots,x_k,x_{k+1},\ldots)=(0,0,\ldots,0,x_k+1,x_{k+1},\ldots)$
where $k$ is the least entry such that $x_k<n_k-1$, and if there is no such $k$ then it produces $(0,0,\ldots)$. Such $(X,T)$ is called an {\em adding machine} or {\em odometer}.

Let $(X,T)$ be t.d.s. A point $x$ is {\em regularly recurrent} if for any open neighborhood $U$ of $x$, there is some $l_U\in \N$ such that $l_U\mathbb{Z}\subset \mathbf{N}(x,U)$. By definition, a regularly recurrent piont is a minimal point. A t.d.s $(X,T)$ is the orbit closure of a regularly recurrent point if and only if  it is an almost one-to-one extension of an adding machine (see \cite{HY2005} or \cite{D}).

If ${\bf \omega}\in \Sigma_k$ is a regularly recurrent point, then ${\bf \omega}$ is called a {\em Toeplitz sequence}, and the shift orbit closure $(X_{\bf \omega},\sigma)$ of a Toeplitz sequence ${\bf \omega}$ is called a {\em Toeplitz system}. Toeplitz systems can be characterized up
to topological conjugacy as topological dynamical systems with the following three properties:
(1) minimal,
(2) almost 1-1 extensions of adding machines,
(3) symbolic.

\subsection{Main result}

It is well-known that for any t.d.s. $(X,T)$ with zero entropy, there exists an extension $(X',\sigma)$ which is a subshift system with zero entropy \cite[Theorem 8.6]{BD2004}. Therefore to prove Sarnak conjecture, one need only verify it for subshift systems with zero entropy. That is,  Sarnak conjecture holds  if for each integer $k\ge 2$ and for any sequence ${\bf a}=(a_n)_{n\in \Z}\in \Sigma_k= \{1,\ldots,k\}^\Z$ with zero entropy,
$$ \lim_{N\rightarrow \infty}\frac{1}{N}\sum_{n=1}^N\mu(n) a_n=0.$$


One motivation of the paper is to find better sequences than the general ones $(a_n)_{n\in \Z}$. The following is the main result of the paper.

\begin{thm}\label{main theorem}
Let ${\bf a}=(a_n)_{n\in \mathbb{Z}}\in \Sigma_k=\{1,\ldots,k\}^{\mathbb{Z}}$. For any $\epsilon>0$, there exists a Toeplitz sequence ${\bf b}=(b_n)_{n\in \mathbb{Z}}\in \Sigma_k$ such that
\begin{enumerate}
  \item $h(X_{\bf b},\sigma)\leq 2h(X_{\bf a},\sigma)$;
  \item $\displaystyle \lim_{N\rightarrow\infty}\frac{1}{2N+1}\sum_{n=-N}^N|a_n-b_n|<\epsilon$.
\end{enumerate}
\end{thm}

By Theorem \ref{main theorem}, one has the following corollary immediately.

\begin{cor}
If the M\"{o}bius function is disjoint from any Toeplitz sequence with zero entropy, then the Sarnak conjecture holds.
\end{cor}

See \cite{DK} for the progress on Sarnak conjecture for Toeplitz systems.

\section{Proof of Theorem \ref{main theorem}}

In this section, we prove the main result Theorem \ref{main theorem}. In the whole section $k\ge 2$ is a fixed integer. First we need some notations.
For $n\in\N$ and words $A_1,\ldots, A_n$, we denote by $A_1\ldots A_n$ the concatenation of $A_1,\ldots, A_n$. Let $(X,\sigma)$ be a subshift and $A$ be a word appeared in $X$. For $s\in \Z$, set
$${_s}[A]_{s+|A|-1}={_s}[A]_{s+|A|-1}^X=\{{\bf x}\in X: x_s x_{s+1}\cdots x_{(s+|A|-1)}=A\},$$
which is called a block of $X$. All blocks of $X$ forms a clopen base of $X$.

\medskip

Let ${\bf a}=(a_n)_{n\in \mathbb{Z}}\in \Sigma_k=\{1,\ldots,k\}^{\mathbb{Z}}$. For any $\epsilon>0$, we will construct a Toeplitz sequence ${\bf b}=(b_n)_{n\in \mathbb{Z}}\in \Sigma_k$ such that
\begin{enumerate}
  \item[(i)] $h(X_{\bf b},\sigma)\leq 2h(X_{\bf a},\sigma)$;
  \item[(ii)] $\displaystyle \lim_{N\rightarrow\infty}\frac{1}{2N+1}\sum_{n=-N}^N|a_n-b_n|<\epsilon$.
\end{enumerate}

To this aim, we will construct a sequence $\{{\bf a}^{(n)}\}_{n=1}^\infty\subset \Sigma_k$, and the sequence ${\bf b}$ is the limit of $\{{\bf a}^{(n)}\}_{n=1}^\infty$, i.e. ${\bf b}= \lim_{n\to\infty} {\bf a}^{(n)}$.

\medskip

Now we construct the  sequence $\{{\bf a}^{(n)}\}_{n=1}^\infty\subset \Sigma_k$ inductively.

\noindent \textbf{Step $1$}: Let $\ep_1<\frac{1}{2}\ep$. Let $l_1\in \mathbb{N}$ such that $\frac{k}{l_1}\leq \ep_1$. Assume that $a_0=j_0$. Let ${\bf a^{(1)}}=\{a_n^{(1)}\}_{n\in \Z}$ be the sequence defined as follows:
$$a_n^{(1)}=\left\{
              \begin{array}{ll}
                j_0, & \hbox{if $n\in l_1\Z$;} \\
                a_n, & \hbox{otherwise.}
              \end{array}
            \right.
$$

Let $[j_0]={_0}[j_0]_0=\{{\bf x}\in X_{{\bf a}^{(1)}}: x_0=j_0\}$.
Then we have the following properties:
\begin{enumerate}
\item[$(I)_1$:] $l_1\Z\subset N({\bf a}^{(1)},[j_0])$;
\item[$(II)_1$:] $\displaystyle \frac{1}{2N+1}\sum_{n=-N}^N|a_n^{(1)}-a_n|\leq \frac{k}{l_1}<\ep_1$ for all $N\geq 1$.
\end{enumerate}

\medskip

Let $\ep_0=\ep$ and $l_0=1$.
For $M\geq 2$, assume that we have ${\bf a}^{(1)}, \ldots, {\bf a}^{(M-1)}\in \Sigma_k$, positive numbers $\ep_1, \ep_2,\ldots, \ep_{M-1}$ and $l_1,l_2,\ldots,l_{M-1}\in \N$ such that for all $j\in \{1,2,\ldots, M-1\}$, $\ep_{j}<\frac 12 \ep_{j-1}$, $l_{j-1}|l_{j}$ with $\frac{2kl_{j-1}}{l_j}\leq \ep_j$ and
\begin{enumerate}
\item[$(I)_j$:] $l_j\Z\subset N({\bf a}^{(j)},[\varpi^{(j)}])$, where $\varpi^{(j)}={\bf a}^{(j-1)}[-l_{j-1}, l_{j-1}-1]=(a^{(j-1)}_{-l_{j-1}},\ldots,a^{(j-1)}_{l_{j-1}-1})$ and $[\varpi^{(j)}]={_{-l_{j-1}}}[\varpi^{(j-1)}]_{l_{j-1}-1}$.
\item[$(II)_j$:] $\displaystyle \frac{1}{2N+1}\sum_{n=-N}^N|a_n^{(j)}-a_n^{(j-1)}|\leq \frac{2kl_{j-1}}{l_j}<\ep_j$ for all $N\geq 1$.
\item[$(III)_j$:] $\varpi^{(j)}=(a^{(j-1)}_{-l_{j-1}},\ldots,a^{(j-1)}_{l_{j-1}-1})
    =(a^{(j)}_{-l_{j-1}},\ldots,a^{(j)}_{l_{j-1}-1})$.
\end{enumerate}

\medskip

Now we give ${\bf a}^{(M)}$.

\noindent \textbf{Step $M$:} Let $0< \ep_M<\frac{1}{2}\ep_{M-1}$. Choose $l_M\in \mathbb{N}$  such that $\frac{2kl_{M-1}}{l_M}\leq \ep_M$ and $l_{M-1}|l_M$ . Let $$\varpi^{(M)}={\bf a}^{(M-1)}[-l_{M-1}, l_{M-1}-1]=(a^{(M-1)}_{-l_{M-1}},\ldots,a^{(M-1)}_{l_{M-1}-1}).$$
Define ${\bf a}^{(M)}\in \Sigma_k$  as follows:
\begin{equation}\label{equation-construction}
   \left\{
              \begin{array}{ll}
                (a_{rl_M-l_{M-1}}^{(M)},\ldots,a_{rl_M+l_{M-1}-1}^{(M)})=\varpi^{(M)}, & \hbox{for all $r\in \mathbb{Z}$;} \\
                a_n^{(M)}=a_n^{(M-1)}, & \hbox{$n\not\in \bigcup_{r\in \Z}[rl_M-l_{M-1}, rl_M+l_{M-1}-1]$.}
              \end{array}
            \right.
\end{equation}

Then we have the following properties:
\begin{enumerate}
\item[$(I)_M$:] $l_M\Z\subset \mathbf{N}({\bf a}^{(M)},[\varpi^{(M)}])$, where $[\varpi^{(M)}]={_{-l_{M-1}}}[{\bf a}^{(M)}]^{X_{\bf b}}_{l_{M-1}-1}$
\item[$(II)_M$:] $\displaystyle \frac{1}{2N+1}\sum_{n=-N}^N|a_n^{(M)}-a_n^{(M-1)}|\leq \frac{2 kl_{M-1}}{l_M}<\ep_M$ for all $N\geq 1$.
 \item[$(III)_M$:] $\varpi^{(M)}=(a^{(M-1)}_{-l_{M-1}},\ldots,a^{(M-1)}_{l_{M-1}-1})=(a^{(M)}_{-l_{M-1}},\ldots,a^{(M)}_{l_{M-1}-1})$.
\end{enumerate}

Thus by induction, we have sequence $\{{\bf a}^{(n)}\}_{n=1}^\infty\subset \Sigma_k$. Note by $(III)_M$, we have that
\begin{equation}\label{equation-IV}
  (a^{(M')}_{-l_{M-1}},\ldots,a^{(M')}_{l_{M-1}-1})
  =(a^{(M)}_{-l_{M-1}},\ldots,a^{(M)}_{l_{M-1}-1}), \ \text{for all} \ M'>M,
\end{equation}
and hence $\lim_{M\rightarrow \infty}{\bf a}^{(M)}$ exists. Let ${\bf b}= (b_n)_{n\in \mathbb{Z}}= \lim_{n\to\infty} {\bf a}^{(n)}$.

\medskip

By Property $(II)_M$, we have that $${\displaystyle\lim_{N\rightarrow\infty}\frac{1}{2N+1}\sum_{n=-N}^N|a_n-b_n|
<\sum_{n=1}^\infty \ep_n<\epsilon}.$$

By (\ref{equation-IV}) and the definition of ${\bf b}$,
we have that
\begin{equation*}
   (b_{rl_M-l_{M-1}},\ldots,b_{rl_M+l_{M-1}-1})=\varpi^{(M)}
  =(a^{(M)}_{-l_{M-1}},\ldots,a^{(M)}_{l_{M-1}-1}), \forall r\in \Z.
\end{equation*}
Hence
$$l_M\Z\subset \mathbf{N}({\bf b},[\varpi^{(M)}]),$$ where  $[\varpi^{(M)}]^{X_{\bf b}}={_{-l_{M-1}}}[{\bf b}]^{X_{\bf b}}_{l_{M-1}-1}$. Since $\lim_{M\to\infty }l_M=\infty$, $\{[\varpi^{(M)}]^{X_{\bf b}}\}_{M=1}^\infty$ is a base of ${\bf b}$ in $X_{\bf b}$ and it follows that ${\bf b}=(b_n)_{n\in \mathbb{Z}}$ is a Toeplitz sequence .

\medskip

It remains to show that $h(X_{\bf b},\sigma)\leq 2h(X_{\bf a}, \sigma)$.

For a sequence ${\bf q}=(q_n)_{n\in \mathbb{Z}}\in \Sigma_k= \{1,\ldots,k\}^{\mathbb{Z}}$, define
\begin{equation}\label{collection of mm word}
W_M({\bf q})=\{(q_{rl_{M}},\ldots,q_{rl_M+l_M-1}): r\in \mathbb{Z}\}\subset \{1,2,\ldots,k\}^{l_M},
\end{equation} i.e. $W_M\big((q_n)\big)$ is the collection of $l_M$-length word which appears in the position $[rl_{M},rl_M+l_M-1]$ for some $r\in \mathbb{Z}$.

 As $l_{M-1}|l_M$, by $(III)_M$ one has that for $M'> M$, $\varpi^{(M')}$ consists of words in $W_{M}({\bf a}^{(M)})$. By (\ref{equation-construction}), one has that
\begin{equation}\label{M=M'}
W_{M}({\bf a}^{(M')})=W_{M}({\bf a}^{(M)}), \ \text{for all }\ M'\ge M.
\end{equation}

Notice that (\ref{equation-construction}) gives a surjective map from $W_{M}({\bf a}^{(M-1)})$ to $W_{M}({\bf a}^{(M)})$. Indeed, for $M'<M$, as $l_{M'}|l_M$, (\ref{equation-construction}) for $M'$ gives a surjective map from $W_{M}({\bf a}^{(M'-1)})$ to $W_{M}({\bf a}^{(M')})$. Therefore
\begin{equation}\label{M>M'}
\# W_{M}({\bf a}^{(M')})\ge \# W_{M}({\bf a}^{(M)}), \ \text{for all } \ M'< M.
\end{equation}

By Equation \eqref{M=M'} and \eqref{M>M'},
$$\# W_{M}({\bf a})\geq \# W_{M}({\bf a}^{(M)})\geq  \# W_{M}({\bf b}), \forall M\in\N.$$

Recall that for a sequence ${\bf q}=(q_n)\in\Sigma_k= \{1,\ldots,k\}^{\Z}$,
$$B_n({\bf q})=\{A: A \text{ is the } n \text{-length word appeared in }{\bf q} \}\subset \{1,2,\ldots,k\}^n,$$ and
\begin{equation*}\label{entropy for (rn)}
h(X_{\bf q},\sigma)=h({\bf q})= \lim_{n\rightarrow \infty } \frac{1}{n}\log\# B_n({\bf q}).
\end{equation*}
Notice that $W_M({\bf a})\subset B_{l_M}({\bf a})$. Therefore
\begin{equation}\label{inequation for words}
\# W_M({\bf b}) \leq \# W_M({\bf a})\leq \# B_{l_M}({\bf a}).
\end{equation}
Notice that for any $A=(w_1,\ldots,w_{l_M}) \in B_{l_M}({\bf b})$, there exists $r\in \Z$ such that $\w$ appears in
$${\bf b}[rl_M, (r+2)l_M-1]=(b_{rl_{M}},b_{rl_{M}+1},\ldots,b_{(r+2)l_M-1}).$$
For any $r\in \Z$, there are at most $l_M+1$ different $l_M$-length words appeared in the $2l_M$-length word $(b_{rl_{M}},b_{rl_{M}+1},\ldots,b_{(r+2)l_M-1}).$
As $(b_{rl_{M}},\ldots,b_{(r+1)l_{M}-1}),(b_{(r+1)l_{M}},\ldots,b_{(r+2)l_{M}-1})\in W_M({\bf b})$, one has that
$$\# B_{l_M}({\bf b}) \leq (l_M+1) \cdot \left(\# W_M({\bf b})\right)^2.$$
Together with \eqref{inequation for words}, one has that
\begin{equation*}\label{conclusion for entropy}
\begin{split}
h({\bf b})&=\lim_{M\to\infty} \frac{1}{l_M} \log \# B_{l_M}({\bf b})\le \lim_{M\to\infty} \frac{1}{l_M} \log[ (l_M+1)\left(\# W_M({\bf b})\right)^2 ] \\ &\leq \lim_{M\rightarrow \infty }\frac{1}{l_M}\log[((l_M+1)\left(\# B_{l_M}({\bf a})\right)^2]\\ &=\lim_{M\to\infty} \frac{\log (l_M+1)}{l_M}+ 2\lim_{M\rightarrow  \infty }\frac{\log\left(\# B_{l_M}({\bf a})\right)}{l_M} =2h({\bf a}).
\end{split}
\end{equation*}
The proof is completed.

\end{document}